\documentclass[journal,twocolumn,letterpaper]{IEEEJERM}
\ifCLASSINFOpdf
\else
\fi
\usepackage{times,amsmath,epsfig}
\usepackage{fancyhdr}
\usepackage{amsmath}
\usepackage{amsfonts}
\usepackage{amssymb}
\usepackage{array}
\usepackage{graphicx}
\usepackage{url}
\usepackage{subfigure}
\usepackage{bm}
\usepackage{breqn}
\usepackage{xcolor}
\usepackage{soul}
\usepackage{amssymb}
\usepackage{stfloats}
\usepackage{subfig}
\usepackage{caption}
\usepackage{subfloat} 
\usepackage[listings]{tcolorbox}
\usepackage{pifont}
\usepackage{url}
\usepackage{graphicx}
\usepackage{color}
\usepackage{placeins}
\usepackage{floatrow}
\usepackage[ansinew]{inputenc}
\usepackage{cancel}

\newfloatcommand{capbtabbox}{table}[][\FBwidth]

\usepackage{blindtext}

 \usepackage{subfig}
 \usepackage{caption}

\usepackage{array, tabu, multirow}
\usepackage[backend=bibtex,style=numeric,sortcites,natbib=true,sorting=none]{biblatex}
\begin{document}

\title{Numerical Modeling for Shoulder Injury Detection Using Microwave Imaging}

\author{Sahar Borzooei,
       Pierre-Henri Tournier,
       Victorita Dolean,
       Christian Pichot,
       Nadine Joachimowicz,
       Helene Roussel,
       Claire Migliaccio
       % <-this % stops a space
}% <-this % stops a space

\markboth{IEEE Journal of Electromagnetics, RF and Microwaves in Medicine and Biology}
{Z. Peng \MakeLowercase{\textit{et al.}}: Portable Electromagnetic Imaging System for Shoulder}

\twocolumn[
\begin{@twocolumnfalse}
  
\maketitle

\begin{abstract}
A portable imaging system for the on-site detection of shoulder injury is necessary to identify its extent and avoid its development to severe condition. Here, firstly a microwave tomography system is introduced using state-of-the-art numerical modeling 
and parallel computing for imaging different tissues in the shoulder.  The results show that the proposed method is capable of accurately detecting and localizing rotator cuff tears of different size. In the next step, an efficient design in terms of computing time and complexity is proposed to detect the variations in the injured model with respect to the healthy model. The method is based on finite element discretization and uses parallel preconditioners from the 
domain decomposition method to accelerate computations. It is implemented using the open source FreeFEM software.
\end{abstract}


\end{@twocolumnfalse}]

% Put footnotes here
{
  \renewcommand{\thefootnote}{}%
  \footnotetext[1]{Sahar Borzooei, Victorira Dolean, Christian Pichot and Claire Migliaccio are with Universit\'e C\^ote d'Azur, Nice, France. (email: sahar.borzooei@etu.univ-cotedazur.fr)}
  \footnotetext[2]{Pierre-Henri Tournier is with Sorbonne Universit\`e, CNRS, Universit\'e Paris Cit\'e, Inria, Laboratoire Jacques-Louis Lions (LJLL), F-75005 Paris, France.}
  % \footnotetext[3]{Claire Migliaccio is with Universit\'e Nice Sophia Antipolis,CNRS, LEAT, nice, France}
  %\footnotetext[4]{Victorita Dolean is with Universit\'e Nice Sophia Antipolis, CNRS, LJAD, Nice, France.
  }
  \footnotetext[3]{H\'el\`ene Roussel is with Sorbonne Universit\'e, CNRS, Laboratoire de G\'enie Electrique et Electronique de Paris, 75252, Paris, France.}
 \footnotetext[4]{Nadine Joachimowicz is with Sorbonne Universit\'e, CNRS, Laboratoire de G\'enie Electrique et Electronique de Paris, 75252, Paris, France, Universit\'e Paris Cit\'e, F-75006 Paris France.}
%}
 
% \fi
\IEEEpeerreviewmaketitle

\section{Introduction}

\begin{figure}[ht]

\centering
\includegraphics[width=0.5\columnwidth]{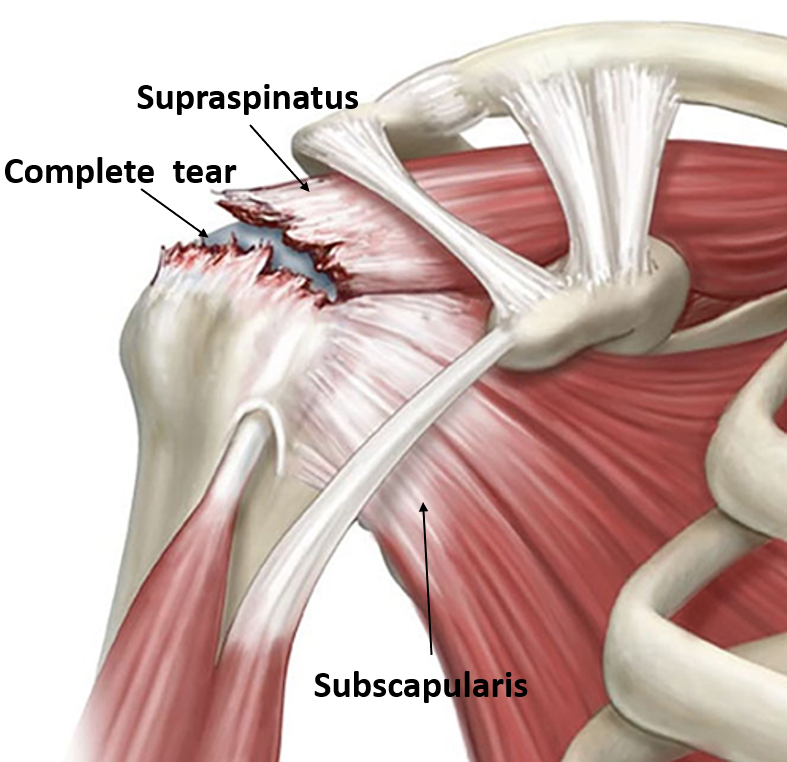}
\caption{Shoulder anatomy~\cite{consult}}
\label{fig_1}
\end{figure}

\IEEEPARstart{R}{otator} cuff tear (RCT) is one of the most common shoulder injuries. It accounts for $70\%$ of
shoulder pain and dysfunction in adults~\cite{Mera}. It can be caused by a variety of factors, including age-related degeneration, overuse and acute injury~\cite{Nyffeler}.
Rotator cuff is a group of muscles and their tendons that work together to stabilize the shoulder, elevate and rotate the arm and keep the head of the humerus securely in the shoulder socket~\cite{Dugas}. 
 Most of the RCTs occur in the supraspinatus tendon~\cite{Omid}. A front view of the shoulder joint and the tendon injury is shown in Figure~\ref{fig_1}~\cite{consult}. 
There are two types of RCTs: partial and complete. A partial tear is when the tendon of the rotator cuff is damaged but not completely severed. A complete tear is when the soft tissue is
detached from the bone. Complete tears can be categorized as small, medium
or large~\cite{Longo}. The severity of RCT does not correlate with the pain experience and it can also be without symptoms, making it challenging to diagnose~\cite{Yamamoto}. 
On the other hand, early detection of a partial RCT can help in preventing its development to full tear and may allow for nonsurgical treatment options~\cite{Abdelwahab}.

Medical imaging plays a critical role in defining the extent of RCTs, which has important implications in clinical decision making and surgical planning~\cite{Konin}. 
Magnetic resonance imaging (MRI), magnetic resonance angiography (MRA) and ultrasound (US) are standard imaging modalities being used to assess the presence and size of RCTs~\cite{Iannotti}. However, these diagnostics techniques are not always accurate in depicting the size and number of involved tendons~\cite{Teefey}. 
MRI as the first-choice imaging modality for the detection of RCTs is costly, bulky and is not suitable for on-site early detection. There is frequent  intra-operative reports claiming to find tears much larger than determined on MRI, or even the lack of a tear~\cite{Xusheng}. MRA is an invasive procedure that requires the injection of a contrast agent~\cite{Yip}. 
US is operator-dependent and may be limited by the patients body habitus~\cite{Adrian}.
As a result, an alternative low-cost, portable and non-invasive method is demanding, specially on-site for competitive athletes such as swimmers~\cite{Stahnke}.
Electromagnetic imaging (EMI) systems have shown promising results in various medical applications~\cite{Kanazawa,Oloumi, Mirbeik, Sultan}. Fully circular tomographic-based EMI systems are designed for different applications, such as brain~\cite{Tournier1, HOPFER}, breast ~\cite{Zellweger} and knee imaging~\cite{Kamel}.
Circular data acquisition effectively helps in increasing the cross-range resolution of the reconstructed images~\cite{Lulu}.\\
To the best of our  knowledge, there is no EMI system to detect shoulder injuries. 
Designing an EMI system for the shoulder is challenging due to following reasons. Firstly, the complex anatomy of the shoulder prevents designing a fully circular phased array for spherical scan and  full data acquisition. As a consequence, the antenna array has to be defined conformal to the shoulder geometry. Another challenge would be the electrically large size of the shoulder along with the heterogeneous nature of the tissues, characterized by high losses (also a characteristic of the matching medium), making it difficult to achieve a high resolution three-dimensional (3D) reconstruction. Shoulder being located in the near-field of the antenna array, an efficient 3D EM-modelling is required to consider coupling effects and near-field interactions between the imaging system and the shoulder. Besides, the variability in the shoulder anatomy among individuals prevents us to consider a particular shoulder structure as a priori knowledge. Finally, the accurate knowledge of the dielectric properties of the shoulder tissues, and more specifically the synovial fluid, remains a challenge.
At the early stage design, when assessing the potential of EMI on a new application, using numerical modeling to help design the system presents several advantages: it allows to accurately model the human body complexity and electrical properties, and to simulate various anatomical scenarios in a flexible way, while reducing cost and time.
In this paper, we propose a feasibility study based on numerical assessments to design a shoulder injury detection system and demonstrate its performance. We make use of an anthropomorphic numerical model of the shoulder and a built-in EM modeling based on the open source code FreeFEM, which offers flexibility in dealing with complex geometries and electrically large structures.
The paper is organized as follows: Section~\ref{sec1} discusses the dielectric properties of the shoulder tissues, while Section~\ref{sec2} describes the choices for a first shoulder imaging system with a dense array of antennas to achieve high resolution. The mathematical framework of the EM modeling is outlined in Section~\ref{sec3}, and numerical solutions are presented in Section~\ref{sec4}. Section~\ref{sec5} investigates the reduction of the number of antennas in the EMI system.

\section{Dielectric properties of the shoulder tissues}
\label{sec1}

\begin{figure}[ht]
\centering
\includegraphics[width=1.0\columnwidth]{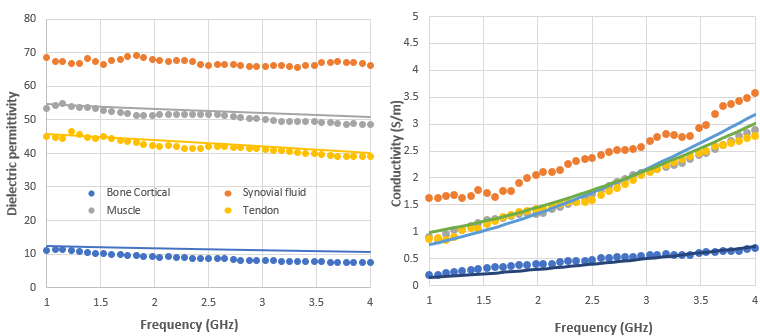}
\caption{Dielectric properties of the tissues}
\label{fig_eps}
\end{figure}

Once the RTCs occurs, synovial fluid (SF) accumulates at the location of the tear~\cite{Papatheodorou, Stone} that leads to a change in the dielectric properties in the shoulder joint~\cite{Kiel}.
It is reported that the volume of SF on aspiration prior to arthroscopic rotator cuff repair correlates with tear size. The mean aspirate volume of partial tears is $ 0.76\pm0.43~ mL$, small tears $ 1.46\pm 1.88~mL~$, medium tears $ 3.04~\pm 2.21 mL~$, and large tears $ 6.60~\pm 3.23~mL$~\cite{Stone}.
 When there is a difference in the dielectric properties between different shoulder tissues and the SF, microwave imaging can detect this contrast. The larger the difference, the easier the detection. Thus, it is essential to measure the dielectric properties of the SF in relation to the other tissues.
In a recent work with the purpose of detecting knee injury~\cite{Kamel}, high values of SF are considered compared to the rest of the tissues; however, to the best of our knowledge, no measurement of dielectric properties of SF is published in this or any other literature. It is reported that low frequencies ranging between $0.5-3~GHz$ are suitable to achieve deeper tissue penetration as well as acceptable imaging quality \cite{nikolova2011}.
In this section, the dielectric properties of SF, muscle, bone cortical and tendons are measured in the frequency range of $1-4~GHz$ that almost include and extend the recommended frequency range. In the case of SF, four patients with knee or shoulder
pain were sampled to obtain the real human SF. We computed the average value of the measured dielectric properties of the four SF samples, which had different temperature and pathology.
The mixtures of the other tissues were made based on the SUPELEC RECIPES website~\cite{SUPELEC}. An optimization code based on Kraszweski's binary law gives the concentration of TritonX-100 and salt required to produce mixtures, whose dielectric properties are close to those given by a Cole-Cole model, for each biological tissue of the shoulder region~\cite{gabriel1996}. Measurements were performed using a homemade coaxial probe connected to a ZVH8 Vector Network Analyzer~\cite{abedi2023}. Fig.~\ref{fig_eps} gathers  all measurements and comparisons with  values from reference websites \footnote{\url{https://itis.swiss/database}},\footnote{\url{http://www.ifac.cnr.it}}. Our approach is validated by the good agreement with the references.

\begin{figure}[ht]
\includegraphics[width=1.0\columnwidth]{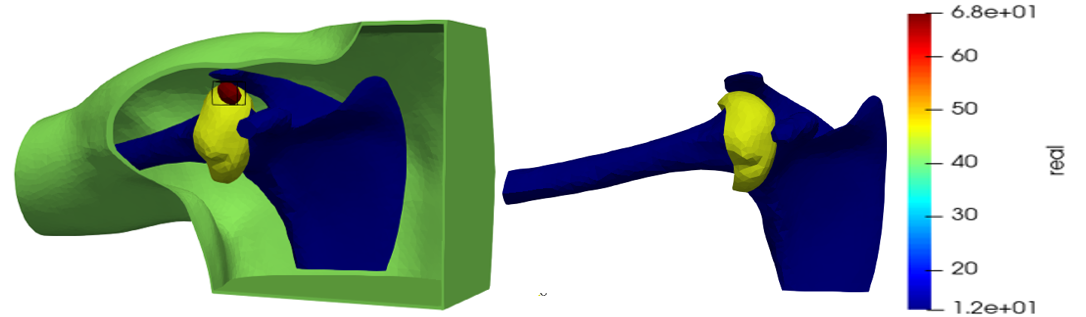}
  \caption{Shoulder Phantoms, real permittivity values}
  \label{fig-2b}
\end{figure}

The relative contrast between SF and other tissues over the frequency range of interest is around 30\%, except for the muscle for which there is a 20\% difference in the real part, making the tear detection feasible but challenging. The anthropomorphic shoulder model is shown in Figure~\ref{fig-2b}. The complex permittivity values of the shoulder tissues at $1~GHz$ are given in Table~\ref{tab_1b} and are used in the simulations. In this Table, the SF value corresponds to our measurements, but the values of other tissues have been already available and are taken from reference websites$^{1},^{2}$. The complex permittivity of the matching medium is chosen equal to that of the muscle.

\begin{table}[ht]
\scalebox{0.9}{
\caption{Complex dielectric properties at $1~GHz$}
		\begin{tabular}{|c|c|c|}
		\hline
		Region &Complex Permittivity& Color in Fig.~\ref{fig-2b}\\
		\hline
		Bone cortical & $12.4 - 2.79j$ & Blue \\
		Tendon & $45.6 - 13.66j$ & Yellow\\
		Muscle& $54.8 - 17.43j$ & Transparent\\
		Skin & $40.9 - 16.17j$ & Green\\
		SF& $68.0 - 29.0j$ & Red\\
		\hline
	\end{tabular}
	\label{tab_1b}
	}
\end{table}

\section{Electromagnetic imaging System}
\label{sec2}
\begin{figure}[ht]
\centering
\includegraphics[width=0.7\columnwidth]{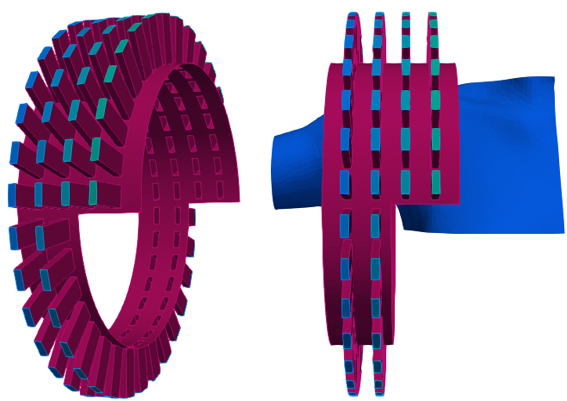}
\caption{EMI system for the shoulder}
\label{fig_2}
\end{figure}

We first design an EMI system with a dense array of antennas  that illuminates the shoulder from different angles. This multi-perspective approach helps in reconstructing a comprehensive and accurate representation of the internal structures~\cite{Liu2019}. Results will be used as reference for the optimized system described in section~\ref{sec5}.
The EMI system consists of 96 ceramic ($\varepsilon_r = 59$) loaded open-ended waveguides arranged on $2$ metallic fully-circular  and $2$ metallic half-circular layers. The two sides of the imaging chamber are open. 
This wearable imaging system with open sides is adapted to the real shoulder structure and designed to surround it partially, as shown in Figure~\ref{fig_2}.
The width and height of the rectangular waveguides are $2.1~cm$ and $0.075~cm$, respectively. Their frequency bandwidth is $0.93-1.85~GHz$.
Here, the operating frequency of $1GHz$ is chosen as a good compromise regarding penetration depth and low specific absorption rate (SAR)~\cite{Sultan1}.
Considering the matching medium as the reference permittivity, the wavelength in this medium is $\lambda = 4.2~cm $. The diameter of the chamber is $7.14\lambda $, the larger length is $3\lambda $ and the distance between antenna layers is $0.47\lambda$. 
The larger size of the modeled
large and partial tears is $0.69\lambda$ and $0.397\lambda$, respectively. Note
that we can expect a better resolution than $0.25\lambda$, because EMI
system is operating in the near field \cite{cui2004study}. The space between antennas and the shoulder is filled with a matching medium to overcome air-skin reflections~\cite{Hoi-Shun}.

\section{Mathematical Framework}
\label{sec3}

\subsection{Finite Element Mesh Generation}

\begin{figure}[ht]
\centering
\includegraphics[width=0.8\columnwidth]{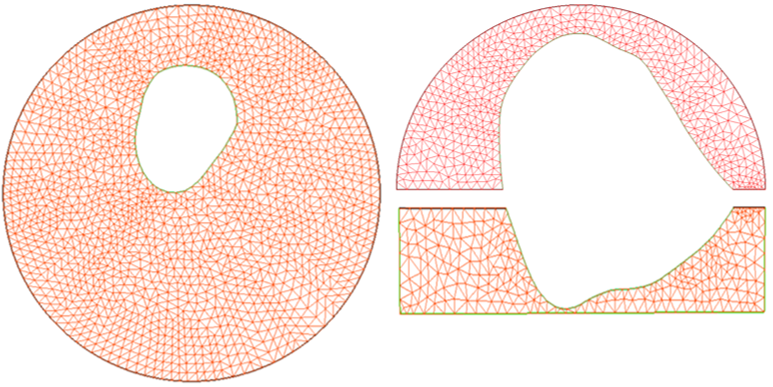}
\caption{The intersecting surfaces of EMI system with shoulder}
\label{fig_mesh}
\end{figure}

Finite element 3D mesh generation of the complete system is a challenging step
due to both the complex geometry of the real body from Computer-Aided Design (CAD) models and the imaging system components. We have chosen to use realistic surface CAD models for shoulder profile and bones including humerus and scapula from a library of 3D models related to the anatomy \cite{Plasticboy}. 
A simple model of rotator cuff tendons was then built surrounding the shoulder joint. The skin is considered with a thickness of $2~mm$ surrounding the muscle geometry.   
The injury is modeled as an ellipsoid in the approximate location of the tear in the supraspinatus tendon, with two different size configurations: small tear with volume of $1 mL$, and large tear with volume of $5 mL$.
Note that in the healthy case, the ellipsoid is replaced with muscle. The regions corresponding to the different tissues are visible in Figure \ref{fig-2b}. The remaining area (excluding bone, injury, tendon, skin) corresponds to muscle. The conformal 3D mesh of the complete system must satisfy non-overlapping elements and nodal-match conditions. At the interfaces between domains, special attention should be given to ensure that the mesh is well-aligned, conformal and continuous~\cite{KRAMER}.
Taking into account that the surface meshes of some tissues intersect the imaging system boundaries, the mesh generation procedure is done by extracting and rebuilding boundary curves corresponding to each intersection for all intersecting surfaces, see Figure \ref{fig_mesh}. In this work, we use the open-source finite element software FreeFEM, which allows to generate 3D meshes and solve partial differential equations (PDEs) with the finite element method (FEM)~\cite{Hecht}. FreeFEM is interfaced with the 
MMG remeshing library [39, 40], which makes it possible to generate adapted tetrahedral meshes of complex geometric configurations, such as the one considered in this paper.

\subsection{Forward Modeling}
\label{secBforward}

The 3D domain ($\Omega$) includes the imaging system and the shoulder as a heterogeneous dissipative  non-magnetic medium of complex permittivity $\varepsilon_{r}= 
(\varepsilon_{r}^{'}- \frac{\sigma j}{\omega \varepsilon_{0}})$, with $\sigma$ the conductivity and $\varepsilon_{r}^{'}$ the electrical permittivity of each tissue, $\varepsilon_{0}$ the electrical permittivity of free space, and $\omega$ the angular frequency.
In the frequency domain, the electric field $\xi(\mathbf{x},\mathbf{t})= \Re(\mathbf{E(x)}e^{\mathrm{i} \omega \mathbf{t}})$ has harmonic dependence on time of angular frequency $\omega$, where $\mathbf{E(x)}$ is its complex amplitude depending on the space variable $\mathbf{x}$. The boundary value problem is defined in equation \ref{eq-m1}.
\begin{equation}
  \begin{cases}
\nabla\times(\nabla\times \mathbf{E})- k^2 \mathbf{E} = 0, & \text{in } \Omega,\\
\mathbf{E}\times \textbf{n} = 0, & \text{on }  \Gamma_m \\
\nabla\times \mathbf{E}\times \textbf{n}+ \mathrm{i}\beta \textbf{n} \times (\mathbf{E}\times \textbf{n}) = g & \text{on }  \Gamma_{t} \\
\nabla\times \mathbf{E}\times \textbf{n}+ \mathrm{i}\beta\textbf{n} \times (\mathbf{E}\times \textbf{n}) = 0 & \text{on }  \Gamma_{r}\\
\nabla\times \mathbf{E}\times \textbf{n}+ \mathrm{i}k \textbf{n} \times (\mathbf{E}\times \textbf{n}) = 0 & \text{on } \Gamma_{o},
  \end{cases}
  \label{eq-m1}
\end{equation}

where $\mathbf{E(x)}$ is the solution of equation~(\ref{eq-m1}) for each transmitting antenna and $k = \omega \sqrt{\varepsilon_{r} \varepsilon_{0}\mu_0}$ is the complex wavenumber of the inhomogeneous medium at each point of the 3D space, with $\mu_{0}$ the permeability of free space. 
$ \beta $ is the propagation wavenumber along the waveguide which corresponds to the propagation of the dominant transverse electric mode $ TE_{10} $. The excitation term is defined as $g= 2\mathrm{i}\beta \mathbf{E}^{TE_{10}}$, imposing an incident wave corresponding to the excitation of the dominant mode of the transmitting waveguide.
$\Gamma_t$ is the port of the transmitting waveguide, $\Gamma_r$ corresponds to ports of receiving waveguides, $\Gamma_m$ represents the metallic surfaces of the waveguides and the walls between them, and $\Gamma_o$ represent the three open sides (right, left and bottom) of the chamber and the boundaries of shoulder profile.   
Let us consider $V = \{\mathbf{v} ~ \in ~ H (\mbox{curl},\Omega),~ \mathbf{v} \times \mathbf{n} =0 ~ \text{on}~ \Gamma_{m}\}$, where $H (\mbox{curl},\Omega)$ is the space of square integrable functions whose curl is also square integrable, and $v$ is the test function. The variational formulation of equation \eqref{eq-m1} is: find $\mathbf{E} ~ \in ~ V $ such that:
\begin{equation} 
\begin{aligned}
\int_{\Omega} [ (\nabla\times \mathbf{E})\cdot (\nabla\times \mathbf{v}) -k^2 \mathbf{E}\cdot \mathbf{v}]
+~\\\shoveleft
 \int_{\Gamma_{t} \cup \Gamma_{r}} \mathrm{i}\beta (\mathbf{E}\times \mathbf{n})\cdot(\mathbf{v} \times \mathbf{n})
 +~\\\shoveleft
 \int_{\Gamma_{o}} \mathrm{i}k (\mathbf{E}\times \mathbf{n})\cdot(\mathbf{v} \times \mathbf{n})& ~=
 \int_{\Gamma_{t}} g \cdot \mathbf{v}, ~
\forall~\mathbf{v} ~ \epsilon ~ V.
\end{aligned}
\label{eq-m2} 
\end{equation}

Through the successive solution of variational problem~\eqref{eq-m2} for each transmitting antenna, we can compute the scattering matrix, a set of complex-valued reflection and transmission coefficients, given in equation~(\ref{eq-s1}):

\begin{equation}
  S_{ij}= \frac{\int_{\Gamma_{r}}\mathbf{E}\cdot\mathbf{E}^{TE_{10}}}{\int_{\Gamma_{r}}|{\mathbf{E}^{TE_{10}}}|^2}.   \label{eq-s1}
  \end{equation}

The computed scattering parameters form a complex matrix of size $96 \times 96$, which is used to produce synthetic data by adding noise to the real and imaginary parts of the coefficients, using a multiplicative white Gaussian noise. The Gaussian noise is applied independently to the real and imaginary parts of each $S_{ij}$ coefficient as independent random variables drawn from the standard normal distribution $N (0,1)$. In this work, we have corrupted the data $S_{ij}$ with $23~dB$ noise.

\subsection{Domain decomposition}
The finite element discretization of our problem leads to a large ill-conditioned linear system $A\mathbf{u}=\mathbf{b}$ for each transmitting antenna. 
Domain decomposition methods (DDMs) are efficient tools to solve such large systems in parallel, both in terms of convergence and computing time~\cite{Bootland}. These
approaches rely on the partition of the computational domain into smaller subdomains,
leading to subproblems of smaller size which are manageable by direct solvers~\cite{Borzooei1}. A Krylov iterative solver (GMRES) along with an Optimized Restricted Additive Schwarz (ORAS) preconditioner is chosen to solve our problem. 
The domain decomposition 
preconditioner is implemented in the HPDDM library~\cite{Jolivet}, an open source high-performance unified framework for domain decomposition methods which is interfaced with the FreeFEM software. 

\subsection{Inverse Problem}
Let $\kappa = k^2$ be the unknown complex parameter of the inverse problem in each point of $\Omega$. In this step an optimization problem, including a fit-to-data term and a regularizing term, is defined with the following cost function:

\begin{equation}\label{sahar:2}
J{(\kappa)}= \frac{1}{2} \sum_{j=1}^{N} \sum_{i=1}^{N} \frac{|S_{ij}(\kappa) - S_{ij}^{syn}|^2}{|S_{ij}^{empty}|^2} + \alpha R(\kappa). 
\end{equation}

Where
\begin{itemize}
 \item $S_{ij}^{syn}$ are the scattering coefficients obtained from the forward problem and are referred to as synthetic data in the rest of the paper.
 \item $S_{ij}^{\kappa}$ are the scattering coefficients computed for the unknown $\kappa$ at the current iteration.
 \item $S_{ij}^{empty}$ are the coefficients computed from the simulation when the domain is filled only with the homogeneous matching medium, used for normalization. 

\item $N$ is the number of antennas of the system.

\item the Tikhonov regularization method is chosen  to reduce the noise effect, defined as $ R(\kappa) = \frac{1}{2}\int_{\Omega}{|\nabla{\kappa}|^2 } $.

\item the regularization parameter $\alpha$ is chosen empirically equal to $10^{-6}$ 
to reach a good compromise between denoising as well as achieving suitable image quality with respect to certain properties such as smoothness or sparsity~\cite{Zheng}.
\end{itemize}

\begin{figure}[ht]
\centering
\includegraphics[width=0.7\columnwidth]{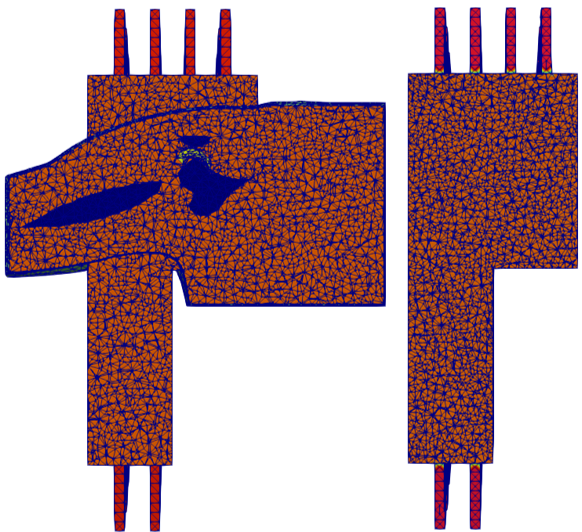}
\caption{Left: Synthetic data, Right: Inverse problem.}
\label{mesh2}
\end{figure}

Minimizing functional~(\ref{sahar:2}) with respect to the parameter $\kappa$ is done by computing its differential and using the adjoint approach in order to simplify its expression \cite{TOURNIER201988}. The gradient of the resulted adjoint functional
is used in the limited-memory Broyden-Fletcher-Goldfarb-Shanno (L-BFGS) optimization algorithm. We refer to~\cite{Tournier1} for a detailed description of the inverse modeling that we have followed in this work. 
Note that we avoid inverse crime by adding noise to the
synthetic data (explained in Section~\ref{secBforward}) and not using a priori knowledge of the body structure. Eliminating a priori knowledge is done by using a mesh that includes the geometry and structure of the body for generating synthetic data (Figure~\ref{mesh2}, left) but defining a different homogeneous mesh that is limited to the imaging chamber (Figure~\ref{mesh2}, right) for the inverse problem.
 Using Nedelec edge finite elements (FE)  of first polynomial degree ($ r=1$) for the FE discretization results in 510531 unknowns for the direct mesh and $357535$ unknowns for the inverse mesh. 
 The spatial resolution for both generated meshes can be defined as the number of points per wavelength $n_{\lambda}$. Here we choose $n_{\lambda} = \frac{\lambda}{6}$, where $\lambda$ corresponds to the wave propagation in the matching medium/muscle tissue. 
 
\section{Numerical Results}
\label{sec4}

\begin{figure}[ht]
\centering
\includegraphics[width=1.0\columnwidth]{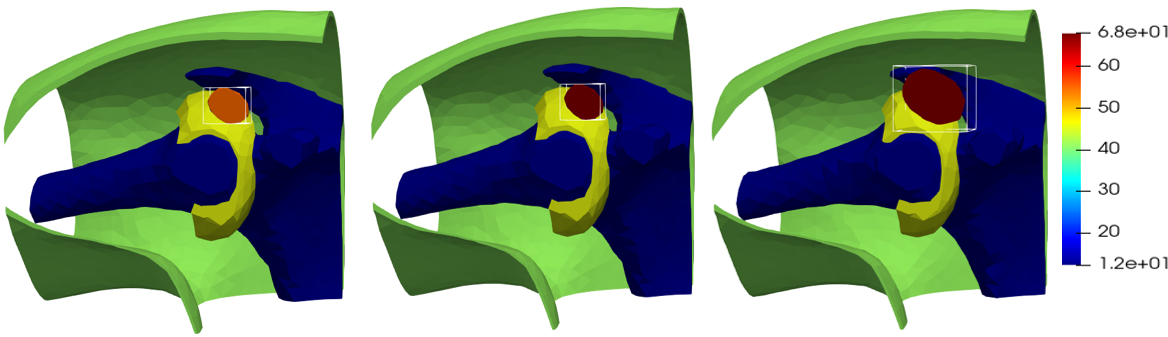}
\includegraphics[width=1.0\columnwidth]{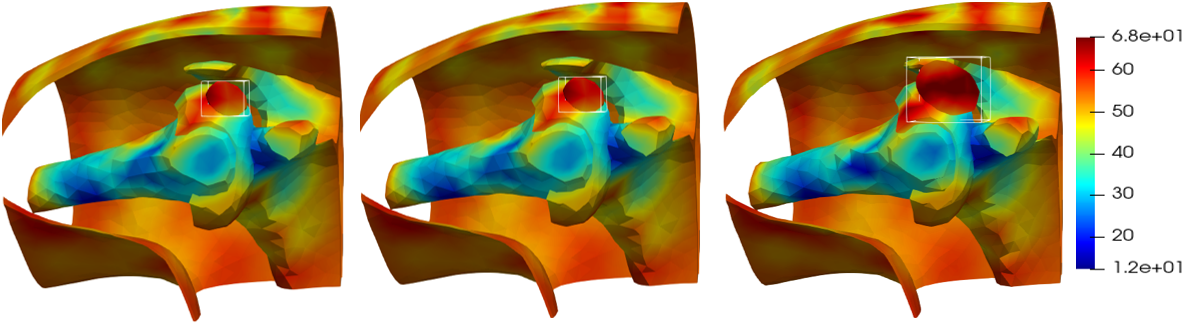}
\caption{Cross-section of the real part of the  permittivity. Top: exact permittivity, bottom: reconstructed permittivity. From left to right: Healthy, partial tear and large tear.}
\label{S-real}
\end{figure}

\begin{figure}[ht]
\centering
\includegraphics[width=1.0\columnwidth]{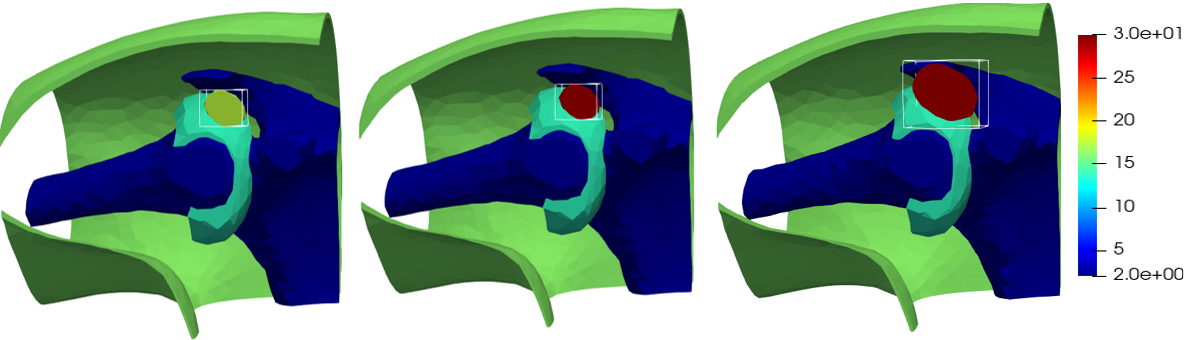}
\includegraphics[width=1.0\columnwidth]{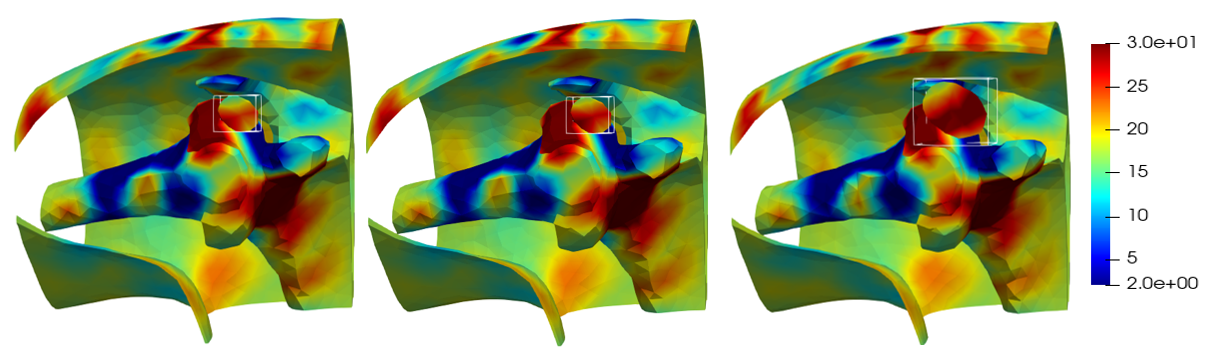}
\caption{Cross-section of the imaginary part of the  permittivity. Top: exact permittivity, bottom: reconstructed permittivity. From left to right: Healthy, partial tear and large tear.}
\label{S-imag}
\end{figure}

\begin{figure}[!htb]
\centering
\includegraphics[width=0.9\columnwidth]{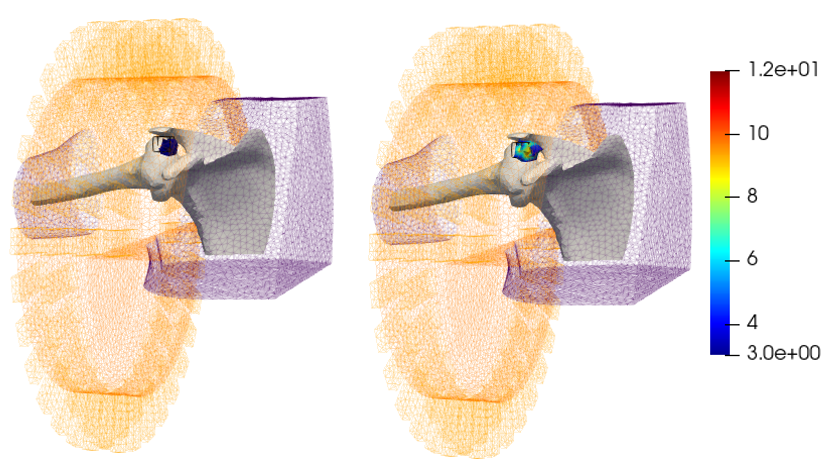}
\caption{The differential image for partial and large tear.}
\label{S-sub}
\end{figure}

Results are obtained on the Universit\'e C\^ote d'Azur's High-Performance Computing (HPC) center. In this HPC center, cluster is composed of 48 CPU computing nodes, including 32 nodes with Dual Intel Xeon Gold processor, providing 40 cores per node and 192 GB of memory and 16 nodes with 2 AMD Epyc processors, providing 32 cores per node and 256 GB of memory. The simulations presented in this paper were carried out using $480$ cores.
%The reconstruction of permittivity is done in full 3D domain, but here the results are shown on the all regions, excluding muscle and matching medium. 
Each reconstruction starts from an initial guess of homogeneous matching medium. The reconstruction results shown here are obtained after 60 iterations; the residual is decreased by a factor of $10^{-2}$. Subsequent iterations do not provide any further noteworthy decrease.
%The solution is converged in $60$ iteration steps when reaching the convergence criterion of $10^{-2}$ for the value of the cost functional. The subsequent iterations only decrease the cost function with the order of $10^{-4}$.
Figures~\ref{S-real} and~\ref{S-imag} show the 3D view of the exact (top) and reconstructed (bottom) results for healthy, partially and fully injured shoulder on the regions of interest. Note that in the healthy shoulder case, the ellipsoid is filled in with muscle as shown in Figs~\ref{S-real}, ~\ref{S-imag} for comparison purpose. The rest of the muscle is not shown for the sake of visibility. The tear is visible in reconstructed real and imaginary parts. For a better view, we also plot the difference between reconstructed permittivity of healthy case and those for the partially and largely injured cases, known as differential images~\cite{Scapaticci}. They are shown for the real part in Figure~\ref{S-sub}, proving that the inverse algorithm succeeds in detecting the injury, with respect to size and location, even for the partial tear. Results for the imaginary part are similar, but are not presented here for the sake of brevity.

\begin{equation}\label{sahar:4}
err_{absolute} = |{\varepsilon_{r}^{reconstructed} - \varepsilon_{r}^{exact}}|
\end{equation} 

\begin{figure}[!htb]
\centering
\includegraphics[width=1.0\columnwidth]{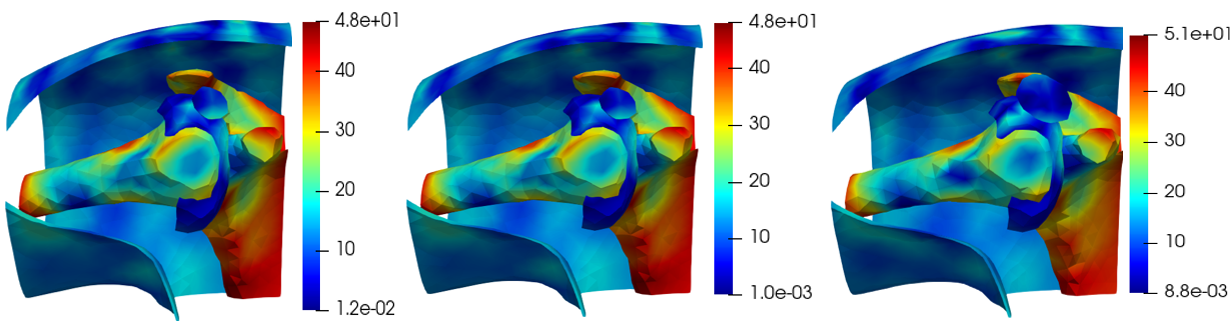}
\caption{The absolute error on real part,  from left to right: Healthy shoulder, partial and large tear}
\label{S-err}
\end{figure}

\begin{figure}[!htb]
\centering
\includegraphics[width=1.0\columnwidth]{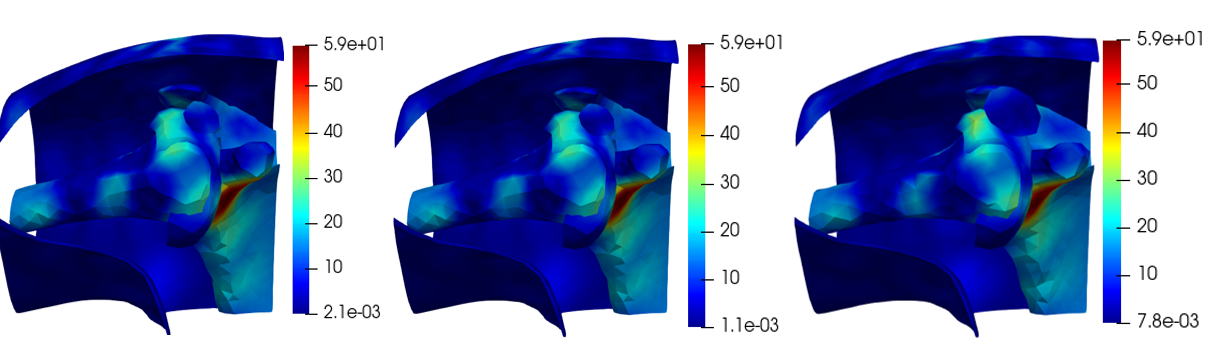}
\caption{The absolute error on imaginary part, from left to right: Healthy shoulder, partial and large tear}
\label{S-err2}
\end{figure}
The reconstruction is further assessed by computing the absolute error on the real and imaginary parts  according to equation~(\ref{sahar:4}) applied on each pixel of the reconstructed domain.  Results are shown in Figure~\ref{S-err} and~\ref{S-err2}.
Error distribution is not uniform across tissues. Highest values are on the edges of the bone area, due to the higher contrast between the complex permittivity of bone and muscle. However, inside the bone in the areas that are not in the vicinity of the edges, the error decreases as expected, due to better performance of the inversion algorithm in a homogeneous region compared to the interfaces. The minimum values of the absolute error in the injury area are $0.78$ and $0.435$, respectively for the real and imaginary parts of the partial tear. These values for the large tear are $0.0088$ and $0.065$, respectively.
\section{An improved design of EMI system }
\label{sec5}

\begin{figure}[h]
\centering
\includegraphics[width=0.8\columnwidth]{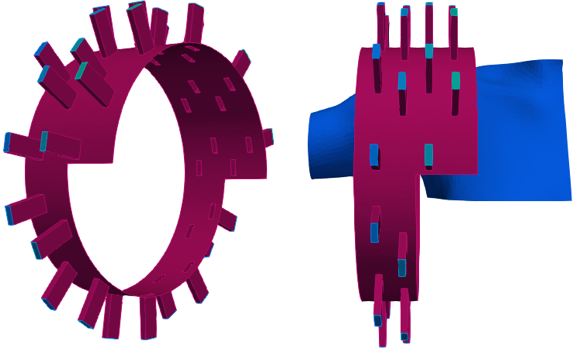}
\caption{Optimal shoulder EMI system}
\label{EMIOPT}
\end{figure}

In this section, an efficient design of the imaging system in terms of reduced number of antennas is proposed. With less transmitting and receiving waveguides, faster computing time as well as cost-effective design can be  achieved.
The system is optimized according to the following guidelines: 
\begin{itemize}
\item The number of antennas is a power of two to ease practically feeding signals into a switching matrix.
\item We introduce a spatial shift between antennas of adjacents rings as shown in Fig.~\ref{EMIOPT}, in order to preserve the spatial diversity. 
\end{itemize}

Different design configurations are assessed based on the value of the $L^2$ norm error, referred to as $err$ according to equation~\eqref{saharerr}, as well as the mean value of the subtraction between the reconstructed injured model and reconstructed healthy model in the injury area, defined in equation~\eqref{saharcont}. Three different configurations are considered, with 16, 32, and 64 antennas ; the reference configuration corresponds to the complete system with $96$ antennas. The values of $err$ for each configuration are computed in the whole domain as well as in the injury area for the partially injured model and are gathered in Table~\ref{table2}. Looking at the $err$ values on both real and imaginary parts, there is no significant increase when reducing the number of antennas by 64 ($N=32$). However, further removing $16$ antennas and going from $N=32$ to $N=16$ antennas leads to a higher increase in $err$. The trend is confirmed by looking at the contrast values $ctr$ which are shown in Table~\ref{table3} for both partial and large tear and with different levels of noise. In the partial tear case, the value of $ctr$ remains positive notable for $N=32$ for all noise levels, and we expect to be able to detect the injury. In contrast, $ctr$ is either negligible or negative for $N=16$, resulting in poor differential images. A first conclusion is that the best design compromise appears to be the $N=32$ configuration, which still allows detection in the most difficult partial tear case. This was further confirmed by extracting and comparing the differential images for each configuration, but due to space limitation the images are not shown in the paper.
The optimal $N=32$ antennas configuration is depicted in Figure~\ref{EMIOPT} and has $11$ and $10$ antennas in both full layers and $6$ and $5$ antennas on both half layers, from left to right. Note that the computing time decreases to about half of the required time for the full $N=96$ configuration.
 
\begin{equation}\label{saharerr}
err = \frac{||\Re(\varepsilon_{r}^{reconstructed}) - \Re(\varepsilon_{r}^{exact})||_2}{||\Re(\varepsilon_{r}^{exact})||_2}
\end{equation} 

\begin{equation}\label{saharcont}
ctr = \varepsilon_{r}^{injured}-\varepsilon_{r}^{healthy}
\end{equation} 

\begin{table}[ht]
\caption{$err$ in the injury area and total domain for the partial tear case, and computing time.}
\scalebox{0.85}{
\begin{tabular}{|c|c|c|c|c|c|c|c|c|}
\hline
\multirow{2}{*}{N}&\multicolumn{2}{c|}{Total}  &\multicolumn{2}{c|}{Injury} & time\\
				& real&imaginary& real&imaginary & (min:s)\\
\hline
		96 &$8.8\%$ & $21.4\%$&$14\%$&$24.4\%$ & $20:09$
 \\
 \hline
 64 &$9.0\%$&$21.2\%$&$10.4\%$&$19.2\%$ & $16:24$
 \\
\hline
		32 &$9.4\%$&$20.6\%$&$16.8\%$&$33.3\%$ & $12:02$
 \\
\hline
		16 &$10.8\%$&$17.2\%$&$19.3\%$&$45.9\%$ &$11:27$
 \\
\hline
	\end{tabular}}
	\label{table2}
\end{table}

\begin{table}[ht]
\caption{$ctr$ in the injury area for partial and large tears with different noise levels.}
\scalebox{0.8}{
\begin{tabular}{c|c|c|c|c|c|c|c|c|c|c|c|}
\cline{2-8}
&\multirow{2}{*}{N}&\multicolumn{2}{c|}{Noise of 23 dB}  &\multicolumn{2}{c|}{Noise of 15 dB}&\multicolumn{2}{c|}{Noise of 10 dB}\\
				&& real&imaginary& real&imaginary& real&imaginary\\
\hline
		 \multirow{4}{*}{\rotatebox{90}{Partial}}&96 &$1.803 $ & $ 1.910$&$ 1.934$&$ 1.767$ & $1.973
 $& $1.696$
 \\
		&64 &$1.605 $&$ 1.673 $&$1.664  $&$ 1.654 $ & $2.233
$ & $ 1.503$
 \\
&32 &$0.890$&$0.690$&$0.843$&$0.669$ &$0.599$&$0.681$
 \\

&16 &$0.133$&$0.1918$&$ -0.136 $&$ 0.0193
$ & $0.3617
 $ & $0.256
 $
 \\
 \hline
 \hline
  \multirow{4}{*}{\rotatebox{90}{Large}}&$96$ & $5.254$& $5.101$& $5.182$ &$ 5.230$ &$5.412$& $5.033$

 \\
&$64$ &$5.053$ &$4.312$ & $5.145$& $ 4.310$& $5.396$
 & $4.244$
 \\
 &$32$ &$2.322$&$2.241$& $2.467$ & $2.191$& $2.524$ & $2.204$
 
 \\
 &$16$& $0.750$ & $0.711$& $0.692$ & $0.681$ & $0.787$ & $0.780$
 \\
 \hline
	\end{tabular}}
	\label{table3}
\end{table}

\section{Conclusion}

We have conducted the first numerical study on the feasibility to image the shoulder joint with an EMI system. The first part of the study consisted in the investigation of the dielectric properties of the shoulder tissues, in particular the synovial fluid which was measured for the first time. Secondly, we have shown that it is possible to reconstruct 3D images of the shoulder joint based on an anthropomorphic numerical model of the shoulder with an imaging system composed of a dense array of $96$ antennas. The reconstruction takes $20$ minutes and $9$ seconds on $480$ computing cores and shows great promise for rapid diagnosis or medical monitoring. Finally, we take advantage of numerical modeling to optimize the number of antennas in the EMI system. We achieve a drastic reduction from $96$ to $32$ antennas. After having proved the relevance of microwave imaging for the detection of shoulder injury, the next step will be to manufacture the imaging system and conduct measurements on built-in phantoms of the different tissues in order to further validate the imaging system proposed in this work.

\section*{Acknowledgment}
This project has received funding from the European Union's Horizon 2020 research and innovation programme under the Marie Sk\l{}odowska-Curie grant agreement No 847581 and is co-funded by the R\'egion Provence-Alpes-C\^ote d'Azur and IDEX $UCA^{JEDI}$and supported by National Research Agency (ANR) under reference number ANR-15-IDEX-01. The authors are grateful to the OPAL infrastructure from Universit\'e C\^ote d'Azur and the Universit\'e C\^ote d'Azur's Center for High-Performance Computing for providing resources and support.
The authors like to thank Dr Eric GIBERT, rheumatologist at the Piti\'e-Salp\^etri\`ere Hospital in Paris, for providing the SF samples for dielectric characterization.

\printbibliography


%\vfill


\end{document}